# Marginal likelihood for parallel series

PETER MCCULLAGH

*Department of Statistics, University of Chicago, Chicago, IL 60637, USA.*
*E-mail: pmcc@galton.uchicago.edu*

Suppose that $k$ series, all having the same autocorrelation function, are observed in parallel at $n$ points in time or space. From a single series of moderate length, the autocorrelation parameter $\beta$ can be estimated with limited accuracy, so we aim to increase the information by formulating a suitable model for the joint distribution of all series. Three Gaussian models of increasing complexity are considered, two of which assume that the series are independent. This paper studies the rate at which the information for $\beta$ accumulates as $k$ increases, possibly even beyond $n$. The profile log likelihood for the model with $k(k+1)/2$ covariance parameters behaves anomalously in two respects. On the one hand, it is a log likelihood, so the derivatives satisfy the Bartlett identities. On the other hand, the Fisher information for $\beta$ increases to a maximum at $k = n/2$, decreasing to zero for $k \geq n$. In any parametric statistical model, one expects the Fisher information to increase with additional data; decreasing Fisher information is an anomaly demanding an explanation.

## 1. Introduction

Let $x_1, \ldots, x_n$ be points in space or time. At each point $x_i$, the $k$-variate response $Y(x_i) = (Y_{i1}, \ldots, Y_{ik})$ is measured. The values are recorded in matrix form $Y = \{Y_{ir}\}$ with one column for each of the $k$ series and one row for each of the $n$ points. Each series is a stationary autoregressive process with autocorrelation parameter $\beta$, and we aim to estimate this parameter as accurately as possible by pooling information from all $k$ series.

Three Gaussian models are considered, all having moments of the form

$$E(Y_{ir}) = 0, \qquad \text{cov}(Y_{ir}, Y_{js}) = \Gamma_{ij} \Sigma_{rs} \tag{1}$$

with autocorrelation function $\Gamma$. The zero-mean assumption is inconsequential and is made for simplicity of notation. It can be replaced by a standard multivariate regression model (Section 3). The three model variants differ only in the assumptions made about the matrix $\Sigma$, which governs the variances and covariances of the $k$ series. These are as follows:

$$\text{Model I: } \Sigma = \sigma^2 I_k, \qquad \text{Model II: } \Sigma = \text{diag}\{\sigma_1^2, \ldots, \sigma_k^2\}, \qquad \text{Model III: } \Sigma \in PD_k,$$

where $PD_k$ is the space of $k \times k$ symmetric positive definite matrices. For each model, we study the profile log likelihood for $\beta$, show that it satisfies the Bartlett identities and study the rate of change of the Fisher information with $k$ and $n$.

Model III aims to accommodate correlations among the series in a simple and natural way, but for $k \geq 2n-1$, the number of parameters exceeds the number of observations. This simple counting argument suggests that we might encounter Neyman–Scott phenomena such as bias, inconsistency or inefficiency in the estimation of $\beta$ (Neyman and Scott (1948)). The failure of profile likelihoods to satisfy the Bartlett identities is the chief explanation for Neyman–Scott phenomena, and the asymptotic bias can often be eliminated by a simple adjustment (Bartlett (1953, 1955), Patterson and Thompson (1971), Cox and Reid (1987), McCullagh and Tibshirani (1990)). The fact that the profile likelihood for $\beta$ in models I–III satisfies the Bartlett identities suggests that Neyman–Scott phenomena should not arise. This intuition is correct for models I and II. However, the marginal likelihood for $\beta$ in model III illustrates a new anomaly for $k > n/2$, namely, that the Fisher information can be increased by deleting one or more series.

Although it is sometimes natural, the separability assumption in (1) is very strong, even for version III. Stein (1999, 2005) is rightly critical of the use of separable covariances for either purely spatial or spatio-temporal processes. However, the product form of the covariance function is extremely convenient and widely used, and there do exist applications in which this assumption is reasonable. It occasionally happens in agricultural field trials that two observations are made on each plot, for example, yield of grain and yield of straw. Although the two yields are certainly correlated, there is good reason to expect that both processes have very similar spatial autocorrelation functions (McCullagh and Clifford (2006)). Mitchell *et al.* (2006) give further references to applications and develop a likelihood-ratio test for separability based on independent replicates of the matrix $Y$. The motivating example for this work arises in a non-spatial context, the estimation of a phylogenetic tree for $n$ species from aligned sequences at multiple homologous loci. Under the model of neutral evolution, the phylogenetic relationship among species is the same at each locus, which implies (1). For further details, see Section 4.

## 2. Profile likelihood

The log likelihood for all three models is

$$\begin{aligned}
l(\Gamma, \Sigma; Y) &= -\frac{1}{2} \log \det(\Gamma \otimes \Sigma) - \frac{1}{2} \text{tr}(Y' \Gamma^{-1} Y \Sigma^{-1}) \\
&= -\frac{k}{2} \log |\Gamma| - \frac{n}{2} \log |\Sigma| - \frac{1}{2} \text{tr}(Y' \Gamma^{-1} Y \Sigma^{-1}),
\end{aligned}$$

using the formula for the determinant of a Kronecker product (Harville (1997), page 350). For fixed $\Gamma$, the log likelihood for model III is maximized at $\hat{\Sigma}_\Gamma = Y'\Gamma^{-1}Y/n$. The log likelihood for model II is maximized at $\mathrm{diag}(\hat{\Sigma}_\Gamma)$ and the log likelihood for model I at $\mathrm{tr}(\hat{\Sigma}_\Gamma)I_k/k$.

The profile log likelihood for $\Gamma$ is

$$l_p(\Gamma;Y) = \begin{cases} -\dfrac{k}{2}\log|\Gamma| - \dfrac{nk}{2}\log\mathrm{tr}(Y'\Gamma^{-1}Y) & \text{(Model I)}, \\ -\dfrac{k}{2}\log|\Gamma| - \dfrac{n}{2}\log|\mathrm{diag}(Y'\Gamma^{-1}Y)| & \text{(Model II)}, \\ -\dfrac{k}{2}\log|\Gamma| - \dfrac{n}{2}\log|Y'\Gamma^{-1}Y| & \text{(Model III)}. \end{cases} \quad (2)$$

The assumption $k \leq n$ is necessary in model III to ensure that the matrix $Y'\Gamma^{-1}Y$ is positive definite with probability one.

The profile log likelihood for model II is a sum over the $k$ series, the contribution of series $r$ being

$$-\frac{1}{2}\log|\Gamma| - \frac{n}{2}\log(Y_r'\Gamma^{-1}Y_r).$$

This is, in fact, the marginal log likelihood based on the standardized statistic $Y_r/\|Y_r\|$, where $Y_r$ is the $r$th column of $Y$ (Bellhouse (1978), Tunnicliffe-Wilson (1989), Cruddas, Reid and Cox (1989)).

For a one-parameter model with derivative matrix $D = d\Gamma/d\beta$, the derivative of the profile log likelihood is

$$2\frac{\partial l_p}{\partial \beta} = \begin{cases} -k\,\mathrm{tr}(WD) + nk\,\mathrm{tr}(Y'AY)/\mathrm{tr}(Y'WY) & \text{(Model I)}, \\ -k\,\mathrm{tr}(WD) + n\displaystyle\sum_{r=1}^{k}(Y_r'AY_r)/(Y_r'WY_r) & \text{(Model II)}, \\ -k\,\mathrm{tr}(WD) + n\,\mathrm{tr}((Y'WY)^{-1}Y'AY) & \text{(Model III)}, \end{cases}$$

where $W = \Gamma^{-1}$ and $A = WDW$. The quadratic form $\mathrm{tr}(Y'WY)$ in model I is distributed as $\sigma^2\chi^2_{nk}$, independently of the ratio $\mathrm{tr}(Y'AY)/\mathrm{tr}(Y'WY)$ (Boos and Hughes-Oliver (1998)). The expected value of the ratio is the ratio of expected values, which is

$$\mathrm{E}\left(\frac{\mathrm{tr}(Y'AY)}{\mathrm{tr}(Y'WY)}\right) = \frac{k\,\mathrm{tr}(A\Gamma)}{nk} = \frac{\mathrm{tr}(WD)}{n}.$$

It follows that the log likelihood derivative for model I has zero expectation. The same argument applied to each series leads to the same conclusion for model II.

The argument for model III is superficially more complicated. For fixed $\Gamma$, the natural quadratic form $Y'WY$ is a complete sufficient statistic for $\Sigma$, with expectation $n\Sigma$. The statistic $\mathrm{tr}((Y'WY)^{-1}Y'AY)$ is invariant under the group $GL(\mathcal{R}^k)$ of linear transformations $Y \mapsto Yg$ acting by right composition. Hence, the distribution does not depend on $\Sigma$. By Basu's theorem (Basu (1955)), every ancillary statistic such as $\mathrm{tr}((Y'WY)^{-1}Y'AY)$ is independent of $Y'WY$. Consequently, if we transform to $Z = W^{1/2}Y$ and condition on the event $Y'WY = Z'Z = I_k$, the columns of $Z$ are orthonormal, the first $k$ columns of a random orthogonal matrix, uniformly distributed with respect to Haar measure on the orthogonal group (Heiberger (1978), Stewart (1980), Diaconis and Shahshahani (1994)). Hence,

$$\begin{aligned}
\mathrm{E}[\mathrm{tr}((Y'WY)^{-1}Y'AY)] &= \mathrm{tr}[\mathrm{E}((Y'WY)^{-1}Y'AY)] \\
&= \mathrm{tr}[\mathrm{E}(Z'\Gamma^{1/2}A\Gamma^{1/2}Z)] \\
&= k\,\mathrm{tr}(\Gamma A)/n \\
&= k\,\mathrm{tr}(WD)/n,
\end{aligned}$$

since $ZZ'$ is a random projection with rank $k \leq n$ and expectation $kI_n/n$. For all three models, the first derivative has zero expectation, so the elimination of $\Sigma$ by maximization has not introduced a bias.

Similar, but more intricate, calculations for random orthogonal matrices described in Appendix A reveal that

$$\mathrm{var}\left(\frac{\partial l_p}{\partial \beta}\right) = -\mathrm{E}\left(\frac{\partial^2 l_p}{\partial \beta^2}\right) = \begin{cases} V\dfrac{k^2}{2(nk+2)} & \text{(I)}, \\ V\dfrac{k}{2(n+2)} & \text{(II)}, \\ V\dfrac{k(n-k)}{2(n-1)(n+2)} & \text{(III)}, \end{cases}$$

where $V = n\,\mathrm{tr}(WDWD) - \mathrm{tr}^2(WD)$. For model III, this formula holds only for $k \leq n$. Thus, the second Bartlett identity is satisfied, and it follows from Appendix B that the Bartlett identities of all orders are satisfied.

For small $k$, the Fisher information increases roughly in proportion to the number of series, all series contributing equally. If, in fact, the series are independent and identically distributed, the efficiency of model II to model I is $(nk+2)/(nk+2k)$, which is fairly high, even for a large number of short series. For example, if $n = 10$, the relative efficiency decreases from 1.0 to $10/12$ as $k \to \infty$. For fixed $k$, the relative efficiency increases with $n$, presumably because the number of nuisance parameters in model II is fixed. It appears from these calculations that the additional flexibility of model II over model I comes at a fairly small cost, so II is likely to be preferred over I in most circumstances.

The most striking anomaly for large $k$ is that the Fisher information for $\beta$ in model III is monotone decreasing for $k > n/2$ and is reduced to zero for $k \geq n$. For a conventional one-parameter model with distributions $f_k(y_1, \ldots, y_k; \beta)$, the Fisher information satisfies

$$\mathrm{FI}_k = \mathrm{FI}_{k-1} + \mathrm{var}\left(\frac{\partial \log f_k(y_k|y_1,\ldots,y_{k-1};\beta)}{\partial \beta}\right) \geq \mathrm{FI}_{k-1},$$

so the Fisher information is necessarily non-decreasing in $k$. It is immaterial whether the components are scalars or vectors. This factorization argument also holds for marginal distributions based on residuals, that is, the REML likelihood for variance components or spatial autocorrelations. It also covers the marginal likelihood for models I and II, and conditional likelihoods of the type used to eliminate nuisance parameters in binary regression models. However, explicit Fisher information calculations for $\beta$ in model III show that this seemingly impregnable argument may fail. The difficulty lies in the fact that the marginal distributions $f_k$ of the maximal invariant in model III cannot be factored: $f_{k-1}$ is not the marginal distribution of $f_k$ under deletion of the last component (see Section 5).

Bearing in mind the stated goal of increasing precision by pooling information from all series, the third formulation is a complete success for small $k$. But it is a spectacular failure for large $k$ because any information about $\beta$ that is present in the first few series remains available even when further series are observed. The marginal likelihood with $k \geq n$ is constant and thus devoid of information, but the marginal likelihood based on any single series or pair of series is informative and the Fisher information is positive. A skeptical reader may consider the case $k = n$, where the matrix $Y$ is invertible with probability one. Direct examination of (2) for model III shows that the term $\det(Y'WY)$ factors and that the log likelihood does not depend on the parameter. These conclusions are independent of the nature of the model for $\Gamma$.

If $k < n$, the log likelihood function for model III may be used for inference about $\beta$, either for computing a point estimate and standard error, for generating confidence intervals or for computing posterior intervals. However, if $k > n/2$, greater precision can be achieved by discarding a random subset of the series and applying the same model to the remainder. This counterintuitive behavior is easily verified by simulation.

## 3. Regression effects

The standard model with zero-mean Gaussian variables is easily extended to include linear models having non-zero mean. The simplest model of this form is the standard Gaussian multivariate regression model,

$$\mathrm{E}(Y) = X\theta, \qquad \mathrm{cov}(Y) = \Gamma \otimes \Sigma, \tag{3}$$

where the model matrix $X$ is of order $n \times p$ with rank $p \leq n$ and $\theta$ is a parameter matrix of order $p \times k$. This model with $k = 2$ occurs in field trials where the response is bivariate, for example, weight of grain and weight of straw on each plot (McCullagh and Clifford (2006)). The log likelihood based on residuals or $X$-contrasts (Patterson and Thompson (1971), Harville (1977)) is

$$\check{l}(\Gamma, \Sigma; Y) = \frac{k}{2} \log \mathrm{Det}(WQ) - \frac{n}{2} \log|\Sigma| - \frac{1}{2} \mathrm{tr}(Y'WQY\Sigma^{-1}),$$

where $Q = I_n - X(X'WX)^{-1}X'W$ has rank $n - p$ and $\mathrm{Det}(\cdot)$ is the product of the non-zero eigenvalues. The profile log likelihood for $\Gamma$ in model III is

$$\check{l}_p(\Gamma; Y) = \frac{k}{2} \log \mathrm{Det}(WQ) - \frac{\max(k, n-p)}{2} \log \mathrm{Det}(Y'WQY).$$

All of the remarks made in the preceding section about the Fisher information hold for the profile residual likelihood with $n$ replaced by $n-p$ and $W$ by $WQ$.

## 4. Application to phylogenetics

The motivating example for this work comes from genetics, where sequence data are observed for $n$ species at $k$ homologous loci. In Kim and Pritchard (2007), $n=5$ and the loci are highly conserved non-coding sequences numbering several thousand. In reality, the value at each locus is a sequence from the genetic alphabet, but we assume here for simplicity that this can be coded in such a way that $Y_{ir}$ is a real number. For locus $r$, the covariance of $Y_{ir}$ with $Y_{jr}$ is $\sigma_{rr}\Gamma_{ij}$, where $\sigma_{rr}$ is the site-specific mutation rate and $\Gamma_{ij}$ is the length of the ancestral tree that is shared by the two species. Under neutral evolution, the genetic distance between species is constant, the same at each locus. Furthermore, the responses at different loci may be correlated due to their proximity on the genomes of one or more species. The natural Gaussian model is (3) with $X = \mathbf{1}$, the constant vector.

Our aim is to estimate the ancestral tree using one of the three variants of (3). The profile log likelihood function for model I is

$$\check{l}(\Gamma; Y) = \frac{k}{2} \log \mathrm{Det}(WQ) - \frac{(n-1)k}{2} \log \mathrm{tr}(Y'WQY)$$
$$= \frac{k}{2} \log \mathrm{Det}(WQ) - \frac{(n-1)k}{2} \log \mathrm{tr}(WQS)$$
$$= \frac{k}{2} \log \mathrm{Det}(WQ) + \frac{(n-1)k}{4} \log \mathrm{tr}(WQD),$$

where $S = YY'$ is the observed inner product matrix and $D_{ij} = S_{ii} + S_{jj} - 2S_{ij}$ is the observed squared distance between species. This expression is the log likelihood function on phylogenetic trees based on the marginal distribution of the squared distance matrix $D$ (McCullagh (2008)).

If we wish to take account of locus-specific mutation rates, version II of the standard model is more appropriate. The profile log likelihood function for this model is

$$\check{l}(\Gamma; Y) = \frac{k}{2} \log \mathrm{Det}(WQ) - \frac{n}{2} \sum_{r=1}^{k} \log(Y_r'WQY_r)$$
$$= \frac{k}{2} \log \mathrm{Det}(WQ) + \frac{n}{4} \sum_{r=1}^{k} \log \mathrm{tr}(WQD_r),$$

which requires locus-specific squared distance matrices $D_r(i,j) = (Y_{ir} - Y_{jr})^2$.

Although formulation III appeared to be appropriate and natural for this application, the model with general $\Sigma$ is a total failure because $k$ is much larger than $n$ and the profile log likelihood is uninformative.

Tractable models intermediate between II and III can be used to take account of correlations and to pool information more efficiently. The technique is illustrated here by the set of Markov matrices, that is, $\Sigma$ is a Green's matrix of the form $a_i b_j$ for $i \leq j$ and $\Sigma^{-1}$ is a symmetric Jacobi, or tri-diagonal, matrix (Karlin (1968), Section 3.3). Let $Y_0 = 0$ and let $Q_r$ be the orthogonal projection in $\mathcal{R}^n$ with kernel $\text{span}(X, Y_{r-1})$ and rank $n - p - 1$ for $r > 1$. Conditional on $Y_1, \ldots, Y_{r-1}$, the residual log likelihood for $\Gamma$ based on $Q_r Y_r / \|Q_r Y_r\|$ is

$$\frac{1}{2} \log \text{Det}(WQ_r) - \frac{\text{rank}(Q_r)}{2} \log(Y_r' W Q_r Y_r).$$

The full log likelihood is the sum of $k$ similar terms, and the derivatives have the same form as those for model II. Since $Q_r$ is a random projection, the matrix $V_r = n \, \text{tr}((WQ_r D)^2) - \text{tr}^2(WQ_r D)$ governing the conditional Fisher information is also random. No closed-form expression is available for the expected value, but symmetry considerations indicate that the total Fisher information is of order $\sum \text{rank}(Q_r)$, directly proportional to the number of series. The marginal likelihood for the series in reverse order is different, but the Fisher information is the same. Neither marginal likelihood coincides exactly with the profile likelihood.

## 5. Marginal likelihood and group orbits

In order to eliminate $\theta$ from the likelihood in the model $Y \sim N(X\theta, \Gamma)$, part of the data is ignored. The residual likelihood function is based on the statistic $LY \sim N(0, L\Gamma L')$, where $L$ is any linear transformation with kernel $\mathcal{X} = \text{span}(X)$. To eliminate scalar constants, we also ignore scalar multiples and base the likelihood function on the reduced statistic $Y/\|Y\|$ or $LY/\|LY\|$. For $k = 1$, this technique gives a marginal log likelihood of

$$\check{l}(\Gamma; Y) = \frac{1}{2} \log \text{Det}(WQ) - \frac{n-p}{2} \log(Y' W Q Y),$$

where $W = \Gamma^{-1}$ and $Q = I - X(X'WX)^{-1} X'W$ has rank $n - p$. Note that $\check{l}(\alpha \Gamma; Y) = \check{L}(\Gamma; Y)$, so the marginal likelihood is constant on scalar multiples of $\Gamma$. Equivalent versions of this marginal likelihood function have been given by Bellhouse (1978, 1990), Cruddas, Reid and Cox (1989) and Tunnicliffe-Wilson (1989).

The marginal log likelihood is based on the maximal invariant under the action of a certain group on the observation space $Y \in \mathcal{R}^n$. The standard residual likelihood associated with the group of translations $Y \mapsto Y + x$ with $x \in \mathcal{X}$ leads to the REML log likelihood

$$\frac{1}{2} \log \text{Det}(WQ) - \frac{1}{2} Y' W Q Y.$$

The maximal invariant can be described in one of two ways, either in terms of $X$-contrasts or in terms of the group orbit which is the coset $y + \mathcal{X}$. When the group is extended to include scalar multiplication, the maximal invariant is reduced and the marginal log likelihood is the function $\breve{l}(\Gamma; Y)$ shown above.

In the multivariate case $Y \sim N(X\theta, \Gamma \otimes \Sigma)$, the regression parameter is eliminated, as above, by considering an arbitrary linear transformation $L: \mathcal{R}^n \to \mathcal{R}^n$ with kernel $\mathcal{X}$ and applying it to each of the columns of $Y$. The kernel is thus $\mathcal{X}^{\oplus k}$, the group orbits are cosets and the multivariate residual log likelihood is

$$\tfrac{k}{2} \log \mathrm{Det}(WQ) - \frac{n}{2} \log |\Sigma| - \tfrac{1}{2} \mathrm{tr}(Y'WQY\Sigma^{-1}).$$

If we now extend the group by linear transformations $Y \mapsto Yg$ with $g \in GL(\mathcal{R}^k)$, the dependence on $\Sigma$ vanishes and the marginal log likelihood is

$$\frac{k}{2} \log \mathrm{Det}(WQ) - \frac{\max(k, n-p)}{2} \log \mathrm{Det}(Y'WQY)$$

(Appendix B). Since this is a log likelihood function, the Bartlett identities are automatically satisfied, as was observed in Section 2.

The preceding remarks help to explain the anomalous behavior of the log likelihood under model III. For $\mathcal{X} = 0$ and $k = 1$, each one-dimensional subspace excluding the origin is a group orbit in $\mathcal{R}^n$, so there are as many orbits as there are points on the projective sphere in $\mathcal{R}^n$. For general $k$, the observation space is $\mathcal{R}^{nk}$, but a typical group orbit has dimension $k^2$, so the maximal invariant has dimension $k(n-k)$, which is the factor governing the rate of increase of the Fisher information. For $k \geq n$, there is one group orbit that has probability one, so the invariant statistic is degenerate and uninformative.

The preceding discussion suggests the following question. The action of the group $GL(\mathcal{R}^k)$ is such that that the maximal invariant has a distribution independent of $\Sigma$. Can the same effect be achieved at less cost by a sub-group? The answer, which is a qualified "yes", is now illustrated by the sub-group $UT_k$ of upper triangular transformations. Taking the series in the order given, the maximal invariant is constructed as follows. For each series $Y_r$, compute the residual after linear regression on both $X$ and $Y_1, \ldots, Y_{r-1}$, ignoring scalar multiples. The contribution to the log likelihood function from the series $Y_r$ is

$$\frac{1}{2} \log \mathrm{Det}(WQ_r) - \frac{\mathrm{rank}(Q_r)}{2} \log(Y_r'WQ_rY_r),$$

where $Q_r$ is the orthogonal projection in $\mathcal{R}^n$ with inner product matrix $W = \Gamma^{-1}$ and null space $\mathrm{span}(X, Y_1, \ldots, Y_{r-1})$. The contribution to the Fisher information is non-negative, but zero for $r \geq n - p$. The total log likelihood based on the maximal invariant under the upper triangular sub-group is thus

$$\sum_{r=1}^{k} \frac{1}{2} \log \mathrm{Det}(WQ_r) - \frac{n-p-r+1}{2} \log(Y_r'WQ_rY_r).$$

The group determines the order in which the series are taken, each order has a different maximal invariant and the log likelihood clearly depends on the order. No closed-form expressions are available for the Fisher information, but, by contrast with the behavior for $GL(\mathcal{R}^k)$, the Fisher information does not decrease with $k$.

## Appendix A: Haar moments

Let $H$ be a random orthogonal matrix uniformly distributed with respect to Haar measure on the orthogonal group of order $n$. The value in row $r$ and column $j$ is denoted by $H_r^j$, so the $(r,j)$ component of $H^2$ is $H_r^i H_i^j$ using the summation convention for repeated indices. By contrast, the $(r,s)$ component of $HH'$ is $H_r^j H_s^j = \delta_{rs}$, where $\delta_{rs}$ is the Kronecker symbol for the identity matrix.

Since $-H$ has the same distribution as $H$, the moments and cumulants of odd order are zero. For $n \geq 2$, the non-zero moments and cumulants up to order four are

$$\mathrm{cov}(H_r^i, H_r^j) = \mathrm{E}(H_r^i H_s^j) = \delta_{rs}\delta^{ij}/n,$$

$$\mathrm{E}(H_r^i H_s^j H_t^k H_u^l) = \frac{(n+1)\delta_{rs}\delta_{tu}\delta^{ij}\delta^{kl}[3] - \delta_{rs}\delta_{tu}\delta^{ik}\delta^{jl}[6]}{n(n-1)(n+2)},$$

$$\mathrm{cum}_4(H_r^i, H_s^j, H_t^k, H_u^l) = \frac{2\delta_{rs}\delta_{tu}\delta^{ij}\delta^{kl}[3] - n\delta_{rs}\delta_{tu}\delta^{ik}\delta^{jl}[6]}{n^2(n-1)(n+2)}.$$

Subscripts on the left-hand sides are in one-to-one alphabetic correspondence with superscripts. For a moment or cumulant of order $k$, the right-hand side is a sum over bi-partitions of $\{1,\ldots,k\}$, that is, ordered pairs of partitions of subscripts and superscripts, all partitions having blocks of size two only. For example, the diagonal bi-partition $(13|24, 13|24)$ appears in alphabetic form as $\delta_{rt}\delta_{su}\delta^{ik}\delta^{jl}$, while $(12|34, 13|24)$ appears as $\delta_{rs}\delta_{tu}\delta^{ik}\delta^{jl}$. The coefficient depends only on the least upper bound of the two partitions. Since there are three partitions of four elements into two blocks of size two, there are nine bi-partitions of $\{1,\ldots,4\}$, the three diagonal elements having one coefficient in the fourth moment and the six off-diagonal elements having a different coefficient. Likewise, there are 15 partitions of six elements into blocks of size two, so the sixth moment is a sum over $15^2 = 225$ bi-partitions. The 15 diagonal pairs have a least upper bound with three blocks, a further 90 pairs have a least upper bound with two blocks and the remaining 120 pairs have a least upper bound with one block. Thus, there are three distinct coefficients in the sum over bi-partitions of $\{1,\ldots,6\}$.

It follows that

$$\mathrm{E}(\mathrm{tr}(H^2)) = \mathrm{E}(H_r^i H_i^r) = \delta_{ir}\delta^{ir}/n = 1,$$

$$\mathrm{E}(\mathrm{tr}^2(H)) = \mathrm{E}(X_r^r X_s^s) = \delta_{rs}\delta^{rs}/n = 1,$$

$$\mathrm{E}(\mathrm{tr}(H^4)) = \mathrm{E}(H_r^u H_s^r H_t^s H_u^t)$$
$$= \delta_{rs}\delta_{tu}((n+1)\delta^{ru}\delta^{st}[3] - \delta^{rt}\delta^{su}[6])/(n(n-1)(n+2))$$
$$= ((n+1)(n^2+2n) - 2n(n+2))/(n(n-1)(n+2)) = 1,$$

$$\mathrm{E}(\mathrm{tr}^2(H^2)) = \mathrm{E}(H_r^s H_s^r H_t^u H_u^t)$$
$$= \delta_{rs}\delta_{tu}((n+1)\delta^{rs}\delta^{tu}[3] - \delta^{rt}\delta^{su}[6])/(n(n-1)(n+2))$$
$$= (3n^2(n+1) - 6n)/(n(n-1)(n+2)) = 3,$$

in agreement with more general formulae for moments of traces given by Diaconis and Shahshahani (1994). Finally, for the variance or covariance of log likelihood derivatives under model III, let $Z$ consist of the first $k$ columns of $H$, so that indices $r, s, \ldots$ run from 1 to $k \leq n$. Then

$$\mathrm{E}\,\mathrm{tr}(Z'AZ) = \mathrm{E}(Z_r^i Z_r^j A_{ij}) = \delta^{rr}\delta^{ij}A_{ij}/n = k\,\mathrm{tr}(A)/n,$$
$$\mathrm{E}(\mathrm{tr}(Z'AZ)\,\mathrm{tr}(Z'BZ)) = \mathrm{E}(Z_r^i Z_r^j A_{ij} Z_s^k Z_s^l B_{kl})$$
$$= A_{ij}B_{kl}\mathrm{E}(Z_r^i Z_r^j Z_s^k Z_s^l)$$
$$= \frac{k(nk+k-2)\,\mathrm{tr}(A)\,\mathrm{tr}(B) + 2k(n-k)\,\mathrm{tr}(AB)}{n(n-1)(n+2)}$$

and

$$\mathrm{cov}(\mathrm{tr}(Z'AZ), \mathrm{tr}(Z'BZ)) = \frac{2k(n-k)(\mathrm{tr}(AB) - \mathrm{tr}(A)\,\mathrm{tr}(B)/n)}{n(n-1)(n+2)}.$$

## Appendix B: Distribution of the maximal invariant

Let $Y$ be a random matrix of order $n \times k$ with density $f(y)\,\mathrm{d}y$ with respect to Lebesgue measure at $y \in \mathcal{R}^{nk}$. In order to calculate the distribution of the maximal invariant under the action of $GL(\mathcal{R}^k)$ by right multiplication, we first observe that the action on the first $k$ components is weakly transitive. For $n = k$, there are many group orbits, but for continuous distributions, there is a single orbit that has probability one. Under standard conditions, the matrix $\tilde{g} = Y^{(k)}$ consisting of the first $k$ rows of $Y$ has full rank, so the group element $\tilde{g}^{-1}$ sends $Y$ to a standard configuration or representative orbit element $Z = Y\tilde{g}^{-1}$ in which the leading $k$ rows are equal to $I_k$.

The Jacobian of the transformation $Y \mapsto (\tilde{g}, Z)$ is equal to $|\tilde{g}|^{n-k}$, so the marginal density of $Z$ is

$$\int_{\mathcal{R}^{k^2}} f(zg)|g|^{n-k}\,\mathrm{d}g.$$

Simplification of this expression for the Gaussian distribution with covariance (1) gives the marginal likelihood function in the form

$$\frac{|\Gamma|^{-k/2}}{|z'\Gamma^{-1}z|^{n/2}} \propto \frac{|\Gamma|^{-k/2}}{|y'\Gamma^{-1}y|^{n/2}}.$$

In other words, the profile log likelihood (2) coincides with the marginal log likelihood based on the maximal invariant.

# Acknowledgement

Support for this research was provided in part by NSF Grant DMS-03-05009.